\documentclass{article}
\usepackage[english]{babel}

\usepackage[letterpaper,top=2cm,bottom=2cm,left=3cm,right=3cm,marginparwidth=1.75cm]{geometry}
\usepackage{mathrsfs}
\usepackage{geometry}
\usepackage{amsmath}
\usepackage{amssymb}
\usepackage{graphicx}
\usepackage[colorlinks=true, allcolors=blue]{hyperref}
\usepackage{amsthm}
\newtheorem{Theorem}{Theorem}[section]
\newtheorem{Corollary}{Corollary}[section]
\newtheorem{Lemma}{Lemma}[section]
\theoremstyle{remark}
\newtheorem{Remark}{Remark}[section]

\allowdisplaybreaks

\title{A sharper  observability for diffusion equations in $\mathbb R^N$}

\author{Yueliang Duan\thanks{Department of Mathematics, Shantou University, Shantou 515063, China;
e-mail: ylduan@stu.edu.cn. This work was supported by the National Natural Science Foundation of China under grant 12201379, the Guangdong Basic and Applied Basic Research Foundation (2026A1515012322) and the Guangdong Provincial Department of Education (2025KCXTD013).} ,
\quad Can Zhang
\thanks{School of Mathematics and Statistics,
Wuhan University, Wuhan 430072, China;
e-mail: canzhang@whu.edu.cn. This work was supported by the
National Natural Science Foundation of China under grant 12422118.}}

\begin{document}
\date{}
\maketitle

\begin{abstract}
We establish that for any small constant $\varepsilon>0$, there exists a nontrivial equidistributed set (with an unbounded complement) such that the observability cost  in the short time $(0,T)$ of diffusion equations over $\mathbb{R}^{N}$ can be improved to a sharper exponential form of $Ce^{\varepsilon/T}$, where $C>0$ is a constant independent of $T$.
\end{abstract}

\medskip


\section{Introduction}
Let $N$ be a positive integer, and let $\mathbb{R}^+=(0,+\infty)$. We consider the following diffusion equation on the whole Euclidean space:
\begin{equation}\label{yu-6-24-1}
\left\{
\begin{array}{lll}
\partial_t u(x,t)-\sum\limits_{j,k=1}^N\partial_j\left(a_{jk}(x)\partial_ku(x,t)\right)=0&\mbox{in}\;\mathbb R^N\times\mathbb R^+,\\
		u(x,0)=u_0(x)&\mbox{in}\;\mathbb R^N,
\end{array}\right.
\end{equation}
where $u_0\in L^2(\mathbb R^N)$, $a_{jk}(x)=a_{kj}(x)$ for all $j,k=1,\cdots,N$, and for all  $x\in\mathbb R^N$. Assume that the diffusion coefficient matrix $A(\cdot)=(a_{jk}(\cdot))_{j,k=1}^N$ is globally Lipschitz continuous, and that it satisfies the uniform ellipticity condition.
That is to say, there is a constant $\lambda\geq1$ such that
\begin{equation*}
|a_{jk}(x)-a_{jk}(y)|\leq \lambda|x-y|\quad\text{and}\quad
\lambda^{-1}|\xi|^2\leq \sum_{j,k=1}^N a_{jk}(x)\xi_j\xi_k \leq \lambda|\xi|^2,
\end{equation*}
for all $j,k=1,\cdots,N$, $x, y\in\mathbb{R}^N$, and for all $\xi=(\xi_1,\xi_2,\cdots,\xi_N)\in\mathbb{R}^N$.

By \cite[Theorem 10.9]{Brezis}, equation \eqref{yu-6-24-1} has a unique solution
$$u\in L^2_{\mathrm{loc}}([0,+\infty); H^1(\mathbb{R}^N))\cap C([0,+\infty); L^2(\mathbb{R}^N)).$$ Moreover, if the initial datum satisfies $u_0\in H^1(\mathbb R^N)$, then the corresponding solution satisfies
$$u\in C([0,+\infty); H^1(\mathbb{R}^N))
\ \text{and}\ \partial_t u\in L^2(\mathbb{R}^+; L^2(\mathbb{R}^N)).$$

There have been many achievements in the study of observability inequalities for solutions of linear parabolic equations in both bounded and unbounded domains $\Omega$. Observability inequalities are closely related to the controllability of the heat equation. Early works include the one-dimensional exact controllability result of Fattorini and Russell in \cite{Fattorini-Russell1971}, and the null controllability result obtained as a singular limit of dissipative wave equations by L$\acute{\mathrm{o}}$pez, Zhang, and Zuazua in \cite{Lopez-Zhang-Zuazua2000}.
When the domain $\Omega$ is bounded, the observability inequalities are obtained by using the Carleman inequality method and the spectral inequality method, as in \cite{Apraiz-Escauriaza22,Apraiz-Escauriaza-Wang-Zhang,Escauriaza-Montaner-Zhang2,Escauriaza-Montaner-Zhang,FI,LEBEAU2,LEBEAU1,LEBEAU3,Phung-Wang-2013,Phung-Wang-Zhang,wang-zhang1}.
In this setting, Miller \cite{Miller2006} established controllability using the control transmutation method, and Ervedoza and Zuazua in \cite{Ervedoza-Zuazua2011} derived sharp observability estimates by transmutation. Using Fursikov-Imanuvilov type Carleman inequalities, B$\acute{\mathrm{a}}$rcena-Petisco in \cite{Barcena-Petisco2020} proved null controllability in pseudo-cylinders with Lipschitz boundary conditions by an internal control.
When the domain $\Omega$ is unbounded, if the observation region $\omega\subset\Omega$ is bounded, the observability inequality may not hold; see, e.g., \cite{MZ,MZb,M05a}. If the observation region is a subset $\omega\subset\Omega$ and $\Omega\backslash\omega$ is bounded, we would like to mention \cite{CMZ}, which establishes sufficient conditions for the observability inequality to hold for heat equations in unbounded domains. The proof there is based on the Carleman estimates. For other related results, the reader may refer to \cite{B,CMV,Gonzalez2007,Z16}.
If the observed region $\omega$ is a thick set, \cite{EV17} and \cite{WangZhangZhang} independently obtained that the observability inequality for the pure heat equation holds true in the whole space. This could be extended to the time-independent parabolic equation associated to the Schr\"odinger operator with analytic coefficients (see \cite{nttv}). The methods employed in these papers are all based on spectral inequalities.

When the coefficients of the parabolic equation are non-analytic, the observability inequality was independently established, through their work on propagation of smallness and spectral estimates, by Duan, Yu, and Zhang in \cite{DYZ}, as well as by Burq and Moyano in \cite{Burq}. It should be emphasized that both studies are built upon the quantitative unique continuation theory for elliptic equations on rough sets, as developed in \cite{LM}. It is worth noting that in all the above results, the observable constant takes the exponential form $Ce^{C/T}$, with a certain constant $C>0$. Surprisingly, by tracking the proof there, even under the condition that the density of thick sets approaches one, it seems that one cannot conclude that the observability constant appearing in the exponent tends to zero.

The motivation of this paper is to establish a sharper observability cost $Ce^{\varepsilon/T}$, which allows $\varepsilon$ to be arbitrarily small by choosing a nontrivial thick set.


In order to state our main results, we first
introduce some notation.
Let $\mathbb N^+=\{1,2,3,\ldots\}$.
For $R>0$ and $x_0\in\mathbb R^N$,
denote by $B_R(x_0)$ the open ball in
$\mathbb R^N$ centered at $x_0$ with
radius $R$.
Similarly, $B_R(x_0,0)$ denotes the open
ball in $\mathbb R^{N+1}$ centered at
$(x_0,0)$ with radius $R$.
We write $X=(x,y)\in\mathbb R^{N+1}$,
where $x\in\mathbb R^N$ and $y\in\mathbb R$.
The usual inner product and norm in
$L^2(B_R(x_0))$ are denoted by
$\langle\cdot,\cdot\rangle_{L^2(B_R(x_0))}$
and $\|\cdot\|_{L^2(B_R(x_0))}$,
respectively. For $r>0$ and $\sigma\in(0,\min\{1,r/2\})$, set
\begin{equation}\label{2026-3a}
 a_\sigma:=\frac{2r+\sigma}{\sqrt N},\qquad
 x_k:=a_\sigma k\quad(k\in\mathbb Z^N),\qquad
 \omega^\sigma:=\bigcup_{k\in\mathbb Z^N}B_r(x_k).
\end{equation}
The balls $B_{r+\sigma}(x_k)$ cover the whole space and the families
$\{B_r(x_k)\}_k$ and $\{B_{2R}(x_k)\}_k$ have uniformly bounded overlap
for every fixed $R>0$.  Moreover, the complement of $\omega^\sigma$
contains open balls arbitrarily far from the origin (see Lemma
\ref{yu-LEMMA-7-10-6} below). $C(\cdot)$ denotes a generic positive constant depending on what are enclosed in the brackets.

Our first main result is the following H\"older-type interpolation inequality, which provides a quantitative unique continuation estimate for solutions of \eqref{yu-6-24-1}.

\begin{Theorem}\label{yu-theorem-7-10-6}
Let $\varepsilon\in(0,1)$ and $r>0$. There exist
constants $C=C(\varepsilon,r,N,\lambda)>0$,
$C'=C'(r,N,\lambda)>0$, and
$\sigma=\sigma(\varepsilon,r,N,\lambda)
\in(0,\min\{1,r/2\})$
such that, for every $T\in(0,1)$ and every
$u_0\in L^2(\mathbb R^N)$, the corresponding
solution $u$ of \eqref{yu-6-24-1} satisfies
\begin{equation*}
\begin{aligned}
\int_{\mathbb R^N}|u(x,T)|^2\,dx
\leq{}&
C e^{\frac{C'\varepsilon}{T}}
\left(
\int_{\omega^\sigma}|u(x,T)|^2\,dx
\right)^{1-\varepsilon}
\left(
\int_{\mathbb R^N}|u_0(x)|^2\,dx
\right)^\varepsilon,
\end{aligned}
\end{equation*}
where $\omega^\sigma$ is the observation set
defined in \eqref{2026-3a}.
\end{Theorem}

As a consequence of Theorem \ref{yu-theorem-7-10-6},  by the telescoping series method, we can obtain the final-state observability inequality in a short time for solutions to the equation \eqref{yu-6-24-1}.
\begin{Theorem}\label{jiudu4}
Let $\varepsilon\in(0,1)$ and $r>0$. Then there exist constants $C=C(\varepsilon,r,N,\lambda)>0$ and $\sigma=\sigma(\varepsilon,r,N,\lambda)\in(0,\min\{1,r/2\})$ such that for any $u_0\in L^2(\mathbb R^N)$ and any $T\in(0,1)$, the corresponding solution $u$ of \eqref{yu-6-24-1} satisfies
\begin{equation*}
\int_{\mathbb R^N}|u(x,T)|^2dx\leq
Ce^{\frac{\varepsilon}{T}}\int_{0}^T\int_{\omega^\sigma}|u(x,t)|^2dxdt,
\end{equation*}
where the observation set $\omega^\sigma$ is as defined in \eqref{2026-3a}.
\end{Theorem}

\begin{Remark}
We do not know how to capture the quantitative dependence of the constants $C$ and $\sigma$ on $\varepsilon$ appearing in the above theorems. This remains an open problem for us.
\end{Remark}

By applying a standard duality argument (see, for instance, \cite{w-08}) to Theorem \ref{jiudu4}, we obtain the following corollary.

\begin{Corollary}
Let $\varepsilon\in(0,1)$ and $r>0$.
There exist constants
$C=C(\varepsilon,r,N,\lambda)>0$ and
$\sigma=\sigma(\varepsilon,r,N,\lambda)
\in(0,\min\{1,r/2\})$
such that, for every $T\in(0,1)$ and every
$y_0\in L^2(\mathbb R^N)$, there exists a control
$f\in L^2(\mathbb R^N\times(0,T))$ supported in
$\omega^\sigma\times(0,T)$ satisfying
\begin{equation*}
\|f\|_{L^2(\mathbb R^N\times(0,T))}
\leq C e^{\frac{\varepsilon}{T}}
\|y_0\|_{L^2(\mathbb R^N)},
\end{equation*}
such that the corresponding solution $y$ of
the controlled system
\begin{equation*}
\begin{cases}
\displaystyle
\partial_t y(x,t)
-\sum_{j,k=1}^N
\partial_j\left(a_{jk}(x)\partial_k y(x,t)\right)
=f(x,t)
& \text{in } \mathbb R^N\times(0,T),\\
y(x,0)=y_0(x)
& \text{in } \mathbb R^N
\end{cases}
\end{equation*}
satisfies $y(x,T)=0$ for a.e. $x\in\mathbb R^N$,
where $\omega^\sigma$ is the observation
region defined in \eqref{2026-3a}.
\end{Corollary}

In addition to the tools developed in \cite{DYZ},
the proof of Theorem \ref{yu-theorem-7-10-6} relies heavily on a new propagation estimate of smallness  for gradients of solutions to elliptic equations established in \cite{NHM}. Precisely, we first prove a local quantitative estimate of unique continuation for solutions to the diffusion equation. Combining with the above local result and the geometry of the equidistributed observation, we then obtain a global quantitative estimate at one time point for solutions of the diffusion equation. We finally apply the telescoping series method to prove Theorem \ref{jiudu4}.

The rest of this paper is organized as follows. In Section \ref{sc2}, we quote three auxiliary lemmas, which are basically taken from \cite{DYZ,NHM}.  Sections \ref{kaodu3} and \ref{finalproof} prove Theorems \ref{yu-theorem-7-10-6} and \ref{jiudu4}, respectively.

\section{Auxiliary lemmas}\label{sc2}
We shall use two lemmas from our earlier work \cite{DYZ}. For the reader's convenience, we will state them below. To this end, fix $0<\tau<T$ and let $\xi\in C^\infty(\mathbb{R}^+;[0,1])$ be a cutoff function satisfying
\begin{equation}\label{yu-6-6-6}
\begin{cases}
	\xi= 1 &\mbox{in}\;[0,\tau],\\
	\xi=0 &\mbox{in}\; [T,+\infty),\\
	|\xi_t|\leq \frac{C}{T-\tau}&\mbox{in}\;(\tau,T),
\end{cases}
\end{equation}
where the positive  constant $C$ is independent  of $\tau$ and $T$. For each $R>0$ and $x_0\in\mathbb R^N$, let $v$ be the solution to
\begin{equation}\label{yu-11-29-4}
\begin{cases}
    \partial_t v(x,t)-\sum\limits_{j,k=1}^N\partial_j\left(a_{jk}(x)\partial_kv(x,t)\right)=0&\mbox{in}\;B_{R}(x_0)\times\mathbb{R}^+,\\
    v=\xi u&\mbox{on}\;\partial B_R(x_0)\times\mathbb{R}^+,\\
    v(\cdot,0)=0&\mbox{in}\; B_R(x_0),
\end{cases}
\end{equation}
where $u$ satisfies \eqref{yu-6-24-1} with the initial value $u_0\in H^1(\mathbb{R}^N)$ and $\xi$ verifies \eqref{yu-6-6-6}. By \cite[Lemma 2.1]{DYZ}, we have the following  exponential decay result.
\begin{Lemma}\label{yu-lemma-6-10-1}
There exist positive constants $C_{1}=C_{1}(\lambda)$ and $C_{2}=C_{2}(\lambda,R,N)$ such that
\begin{equation*}
 	\|v(\cdot,t)\|_{L^2(B_R(x_0))}\leq C_{1}\left(1+\frac{1}{T}\right)^{\frac{1}{2}}
 e^{\frac{C_{1}T}{T-\tau}-C_{2}(t-T)^+}F_{T}(R)\quad\text{for all}\;\; t\in\mathbb{R}^+,
 \end{equation*}
where $(t-T)^+:=\max\{0,t-T\}$, $F_{T}(R):=\sup\limits_{s\in[0,T]}\|u(\cdot,s)\|_{H^1(B_R(x_0))}$.
\end{Lemma}

Define
\begin{equation*}
	\tilde{v}(\cdot,t)=
\begin{cases}
	v(\cdot,t)&\mbox{if}\;\;t\geq 0,\\
	0&\mbox{if}\;\;t<0,
\end{cases}
\end{equation*}
where $v$ is the solution of \eqref{yu-11-29-4}.
By Lemma \ref{yu-lemma-6-10-1}, the Fourier transform of $\tilde{v}(x,\cdot)$ with respect to the time variable is well defined. More precisely, we set
\begin{equation*}
	\hat{v}(x,\mu)=\int_{\mathbb{R}}e^{-i\mu t}\tilde{v}(x,t)dt\quad\text{for}\;\;(x,\mu)\in B_R(x_0)\times\mathbb R
\end{equation*}
is well defined. By \cite[Lemma 2.2]{DYZ}, we moreover have
\begin{Lemma}\label{yu-lemma-6-18-1}
There exists a positive constant $C_{3}=C_{3}(\lambda,N)$ such that for each $\mu\in\mathbb{R}$, the following holds:
\begin{equation*}
	\|\hat{v}(\cdot,\mu)\|_{L^2(B_{\frac{R}{2}}(x_0))}\leq C_{3}
	\left(T+R^{2}\right)e^{\frac{C_{1}T}{T-\tau}-\frac{\sqrt{|\mu|}R}{C_{4}}}F_{T}(R),
\end{equation*}
where $C_{4}:=8\pi e\sqrt{2\sqrt{2+\lambda^{6}}};$ $C_{1}>0$ and $F_{T}(R)$ are given in Lemma \ref{yu-lemma-6-10-1}.
\end{Lemma}

Let $0<r<r'<\infty$ and $h>0$. Set
\begin{equation*}
D_s:=\bigl\{X=(x,x_{N+1})\in B_{r'}(x_0)\times(-h,h):
\operatorname{dist}(X,\overline{B_r(x_0)}\times\{0\})<s\bigr\}.
\end{equation*}
Let $g\in H^1(B_{r'}(x_0)\times(-h,h))$ be a weak solution of
\begin{equation}\label{yu-6-23-9}
\begin{cases} \sum\limits_{j,k=1}^N\partial_j\left(a_{jk}(x)\partial_kg(x,x_{N+1})\right)
+g_{x_{N+1}x_{N+1}}=0&\text{in}\; B_{r'}(x_0)\times(-h,h),\\
	g(x,0)=0&\text{in}\; B_{r'}(x_0),\\
	g_{x_{N+1}}(x,0)=z(x)&\text{in}\; B_{r'}(x_0),
\end{cases}
\end{equation}
where $z\in L^2(B_{r'}(x_0))$.

In order to give the proof of our main result, we need the following auxiliary lemma, which is adapted from \cite[Lemma 2.1]{NHM}.

\begin{Lemma}\label{auxiliary-lemma}
For every $\alpha\in(0,1)$ there exist
$C_5=C_5(r,r',h,\lambda,\alpha)>0$ and
$\delta=\delta(r,r',h,\lambda,\alpha)\in(0,1)$ such that
$D_\delta\Subset B_{r'}(x_0)\times(-h,h)$ and
\begin{equation}\label{2022-2.7}
\|g\|_{H^1(D_\delta)}\leq
C_5\|g\|_{H^1(B_{r'}(x_0)\times(-h,h))}^{\alpha}
\|z\|_{L^2(B_r(x_0))}^{1-\alpha}.
\end{equation}
\end{Lemma}
\begin{proof}
Given $\alpha\in(0,1)$, by applying Lemma 2.1 in \cite{NHM} to $g$, there exist two constants $\delta=\delta(\alpha, r', r,h,\lambda)\in(0,1)$ and $C_{5}=C_{5}(\alpha, r', r, h,\lambda)>0$ such that
\begin{equation*}
	\|g\|_{H^1(D_{\delta})}\leq C_{5}\left(\|g\|_{H^{1/2}(B_{r}(x_0)\times\{0\})}
+\|\partial_{x_{N+1}}g\|_{H^{-1/2}(B_{r}(x_0)\times\{0\})}\right)^{1-\alpha}
\|g\|_{H^1(B_{r'}(x_0)\times(-h,h))}^{\alpha}.
\end{equation*}
This, together with \eqref{yu-6-23-9} and the continuous embedding of \(L^2(B_r(x_0))\) into \(H^{-1/2}(B_r(x_0))\), leads to \eqref{2022-2.7}.
\end{proof}

\begin{Remark}
The idea of proving Lemma \ref{auxiliary-lemma} comes from \cite[Lemma 2.1 and Proposition 1.2]{NHM}.
Indeed, we just slightly modify the presentation of Lemma 2.1 in \cite{NHM}, by changing
$\Omega\times(-1,1)$  and $\omega$ into $B_{r'}(x_0)\times(-h,h)$ and $B_{r}(x_0)$, respectively.
\end{Remark}

\section{Proof of Theorem \ref{yu-theorem-7-10-6}}\label{kaodu3}
In order to give the proof of Theorem \ref{yu-theorem-7-10-6}, we need  a local interpolation inequality (i.e., Lemma \ref{lemma-2A2}). It is worth mentioning that the interpolation inequality  established in Lemma \ref{lemma-2A2} is similar to that obtained in \cite[Lemma 3.1]{DYZ}.
However, the interpolation index $\alpha\in(0,1)$ is arbitrary here,
rather than being some fixed value. This is the key improvement.

Let $\mathcal{A}(\cdot):=-\mbox{div}(A(x)\nabla(\cdot))$
be the Dirichlet realization of the
elliptic operator on $B_R(x_0)$, with domain
\begin{equation*}
D(\mathcal{A}):=H_{0}^1(B_R(x_0))\cap H^2(B_R(x_0)).
\end{equation*}
Using the Poincar\'e inequality, we have that
	there is a constant $C=C(\lambda,N)>0$ such that
 \begin{equation}\label{yu-6-7-9}
 	\langle\mathcal{A}h,h\rangle_{{L}^2(B_R(x_0))}\geq CR^{-2}\|h\|^2_{L^2(B_R(x_0))}
	\;\;\mbox{for each}\;\;h\in D(\mathcal{A}).
 \end{equation}
According to \eqref{yu-6-7-9} and the spectral theorem of operators, there are eigenvalues
$\{\mu_i\}_{i\in\mathbb{N}^+}\subset \mathbb{R}^+$ and eigenfunctions
$\{h_i\}_{i\in\mathbb{N}^+}\subset H_0^1(B_R(x_0))$ such that
 \begin{equation}\label{yu-6-7-10}
 \begin{cases}
 \mathcal{A}h_i=\mu_ih_i\;\mbox{and}\;\|h_i\|_{{L}^2(B_R(x_0))}
 =1&\mbox{for each}\;i\in\mathbb{N}^+,\\
CR^{-2}\leq \mu_1\leq \mu_2\leq \cdots
\leq \mu_i\to+\infty&\mbox{as}\;i\to+\infty.
 \end{cases}
 \end{equation}

\begin{Lemma}\label{lemma-2A2}
Let $T>0$, $r>0$, and $x_0\in\mathbb R^N$.
For every $\alpha\in(0,1)$, there exist constants
$C_6=C_6(r,N,\lambda,\alpha)>0$ and
$\sigma=\sigma(r,N,\lambda,\alpha)
\in(0,\min\{1,r/2\})$
such that, for every $\tau\in(0,T/2]$ and every
initial datum $u_0\in H^1(\mathbb R^N)$,
the corresponding solution $u$ of
\eqref{yu-6-24-1} satisfies
\begin{equation}\label{eq-local-ball}
\begin{aligned}
\|u(\cdot,\tau)\|_{L^2(B_{r+\sigma}(x_0))}
\leq{}&
C_6\left[
(1+T^2)
e^{\frac{2C_1T}{T-\tau}}
+
e^{\frac{2r^2}{\tau}}
\right]^{\frac{1-\alpha}{2}}\\
&\times
\|u(\cdot,\tau)\|_{L^2(B_r(x_0))}^{\alpha}
\left(
\sup_{s\in[0,T]}
\|u(\cdot,s)\|_{H^1(B_R(x_0))}
\right)^{1-\alpha},
\end{aligned}
\end{equation}
where $R=C(N,\lambda)r$, with
$C(N,\lambda)>0$ being a sufficiently large
fixed constant independent of $\alpha$.
The constant $C_1>0$ is given by
Lemma \ref{yu-lemma-6-10-1}.
\end{Lemma}

\begin{proof}
Fix $r'=3r/2$ and $h=r$, and choose
$R=C(N,\lambda)r$ sufficiently large
such that
\begin{equation*}
4r'<R/2,\qquad
R/C_4>\sqrt2 h.
\end{equation*}
Let $u$ be the solution to the equation (\ref{yu-6-24-1})  with the initial value $u_{0}\in H^{1}(\mathbb R^N)$.
Let $u_1$ and $u_2$ be solutions to
\begin{equation*}
\begin{cases}
\partial_tu_{1}+\mathcal{A}u_{1}=0&\mbox{in} \;B_R(x_0)\times(0,2T),\\
u_1=u&\mbox{on}\;\;\partial B_R(x_0)\times(0,2T),\\
u_1(\cdot,0)=0 &\mbox{in}\;\;B_{R}(x_0)
\end{cases}
\end{equation*}	
and
\begin{equation*}
\begin{cases}
\partial_tu_{2}+\mathcal{A}u_{2}=0&\mbox{in}\; B_R(x_0)\times(0,2T),\\
u_2=0&\mbox{on}\;\partial B_R(x_0)\times(0,2T),\\
u_2(\cdot,0)=u_0&\mbox{in}\; B_{R}(x_0),
\end{cases}
\end{equation*}
respectively. By the uniqueness of solutions, we have $u=u_1+u_2$ in $B_R(x_0)\times[0,2T]$.
Fix any  $\tau\in(0,T/2]$. Let $v_1$ be the solution to
\begin{equation}\label{2026-y}
\begin{cases}
\partial_tv_1+\mathcal{A}v_{1}=0&\mbox{in}\;B_{R}(x_0)\times\mathbb{R}^+,\\
v_1=\xi u&\mbox{on}\;\partial B_R(x_0)\times\mathbb{R}^+,\\
v_1(\cdot,0)=0&\mbox{in}\; B_R(x_0),
\end{cases}
\end{equation}
where $\xi$ is given by \eqref{yu-6-6-6}.
It is clear that
$u=v_1+u_2$ in $B_R(x_0)\times[0,\tau]$.
We extend $v_1$ by zero to \(B_R(x_0)\times(-\infty,0)\) and continue to denote the extension by $v_1$.
Define
\begin{equation}\label{2026-y1}
\hat{v}_1(x,z)=\int_{\mathbb{R}}e^{-iz t}v_1(x,t)dt
\ \text{for}\;\;(x,z)\in B_R(x_0)\times\mathbb R.
\end{equation}
According to Lemma \ref{yu-lemma-6-10-1}, $\hat{v}_1$ is well defined.

Let $\kappa_{0}:=\sqrt{2}/C_{4}$ with $C_{4}$ given in Lemma \ref{yu-lemma-6-18-1}.
We define $V=V_1+V_2$ in $B_{R/2}(x_0)\times$ $(-\kappa_{0} R,\kappa_{0} R)$
with
\begin{equation}\label{yu-6-23-6jia}
V_1(x,y)=\frac{1}{2\pi}\int_{\mathbb{R}}e^{iz\tau}\hat{v}_1(x,z) \frac{\sinh(\sqrt{-iz}y)}{\sqrt{-iz}}dz\;\mbox{in}\;B_{R/2}(x_0)\times(-\kappa_{0} R,\kappa_{0} R),
\end{equation}
\begin{equation}\label{yu-7-5-7} V_2(x,y)=\sum_{k=1}^\infty\alpha_ke^{-\mu_k\tau}h_k(x)\frac{\sinh(\sqrt{\mu_k}y)}
	{\sqrt{\mu_k}}\;\mbox{in}\;B_{R/2}(x_0)\times(-\kappa_{0} R,\kappa_{0} R)
\end{equation}
where $\{\mu_k\}_{k=1}^\infty$, $\{h_k\}_{k=1}^{\infty}$ are given by $(\ref{yu-6-7-10})$, and
$\alpha_k=\langle u_2(\cdot,0),h_k\rangle_{L^2(B_R(x_0))}$.
Note from
Lemma \ref{yu-lemma-6-18-1} that $V_1$ is also well defined.
It can be verified that
\begin{equation}\label{2026-xx1}
\begin{cases}
\mbox{div}(A(x)\nabla V(x,y))+V_{yy}(x,y)=0&\mbox{in}\;
B_{\frac{R}{2}}(x_0)\times(-\kappa_{0} R,\kappa_{0} R),\\
V(x,0)=0&\mbox{in}\;B_{\frac{R}{2}}(x_0),\\
V_y(x,0)=u(x,\tau)&\mbox{in}\;B_{\frac{R}{2}}(x_0).
\end{cases}
\end{equation}
First, we verify that $V$ satisfies the elliptic equation. For $V_2$, using \eqref{yu-7-5-7} and \eqref{yu-6-7-10}, we obtain
$$
\operatorname{div}(A\nabla V_2)
= -\sum_{k=1}^{+\infty}\alpha_k e^{-\mu_k\tau}h_k\sqrt{\mu_k}\sinh(\sqrt{\mu_k}y)\;  \mbox{and}\;
(V_2)_{yy}
=\sum_{k=1}^{+\infty}\alpha_k e^{-\mu_k\tau}h_k\sqrt{\mu_k}\sinh(\sqrt{\mu_k}y).
$$
Hence,
\begin{equation}\label{2026-q1}
\operatorname{div}(A\nabla V_2)+(V_2)_{yy}=0.
\end{equation}
For $V_1$, set
\[
w_z(y)=\frac{\sinh(\sqrt{-iz}\,y)}{\sqrt{-iz}},
\]
which satisfies $w_z''=-iz\,w_z$. By \eqref{2026-y} and \eqref{2026-y1}, we have that $\hat v_1$ satisfies $\operatorname{div}(A\nabla\hat v_1)=iz\,\hat v_1$. Combining this with \eqref{yu-6-23-6jia} yields
\begin{equation}\label{2026-q2}
\operatorname{div}(A\nabla V_1)+
(V_1)_{yy}
=0.
\end{equation}
From \eqref{2026-q1} and \eqref{2026-q2}, we obtain
$
\operatorname{div}(A\nabla V)+V_{yy}=0.
$

Second, the condition $V(x,0)=0$ follows immediately since $\sinh(0)=0$ for both $V_1$ and $V_2$.

Third, we compute $V_y(x,0)$. For $V_1$, using
\[
\left.\frac{d}{dy}\frac{\sinh(ay)}{a}\right|_{y=0}=1,
\]
we obtain
\[
(V_1)_y(x,0)=\frac{1}{2\pi}\int_{\mathbb R}e^{i\tau z}\hat v_1(x,z)\,dz=v_1(x,\tau).
\]
For $V_2$,
\[
(V_2)_y(x,0)=\sum_{k=1}^{+\infty}\alpha_k e^{-\mu_k\tau}h_k(x).
\]
Therefore,
\begin{equation}\label{2026-11-13A}
V_y(x,0)=v_1(x,\tau)+\sum_{k=1}^{+\infty}\alpha_k e^{-\mu_k\tau}h_k(x).
\end{equation}
By the decomposition of $u$, we have
\begin{equation}\label{2026-11-13B}
u(x,\tau)=v_1(x,\tau)
+\sum_{k=1}^{\infty}\alpha_k e^{-\mu_k\tau}h_k(x).
\end{equation}
It follows from \eqref{2026-11-13A} and \eqref{2026-11-13B} that
$
V_y(x,0)=u(x,\tau).
$

By Lemma \ref{auxiliary-lemma}, there exist constants $C_5=C_5(r,N,\lambda,\alpha)>0$ and
$
\delta \in(0,
\min\left\{1,
r/8
\right\})
$ such that
\begin{equation}\label{yu-7-4-2}
\begin{aligned}
\|V\|_{H^1(D_\delta)}
\leq{}&
C_5
\|V\|_{H^1(
B_{r'}(x_0)\times(-h,h)
)}^{1-\alpha}
\|u(\cdot,\tau)\|_{L^2(B_r(x_0))}^{\alpha}.
\end{aligned}
\end{equation}
By the interior estimate, there is a constant $C=C(r, N,\lambda)>0$ such that
\begin{equation}\label{yu-7-4-3}
\|V\|_{H^1(B_{r'}(x_0)\times(-h,h))}\leq C \|V\|_{L^2(B_{2r'}(x_0)\times(-2h,2h))}.
\end{equation}
Hence, it follows from  (\ref{yu-7-4-2}) and (\ref{yu-7-4-3}) that
\begin{equation}\label{yu-7-4-1}
\|V\|_{H^1(D_{\delta})}\leq C\|V\|_{L^2(B_{2r'}(x_0)\times(-2h,2h))}^{1-\alpha}
\|u(\cdot,\tau)\|^{\alpha}_{L^2(B_{r}(x_0))}
\end{equation}
with $C=C(r,N,\lambda,\alpha)>0$.
Set $\sigma=\delta/4$ and
$
\rho=\delta/16.
$
In what follows, we establish the estimate
\begin{equation}\label{yu-7-4-1777}
\|u(\tau)\|_{L^2(B_{r+\sigma}(x_0))}
\leq
C(N,\lambda)\delta^{-1/2}
\|V\|_{H^1(D_\delta)},
\end{equation}

For every $x\in B_{r+\sigma}(x_0)$, we have
$
\operatorname{dist}
\bigl((x,0),
\overline{B_r(x_0)}\times\{0\}\bigr)
\leq \sigma.
$
Let $(z,y)\in B_{4\rho}(x,0)$. By the
triangle inequality,
\[
\begin{aligned}
&\operatorname{dist}
\bigl((z,y),
\overline{B_r(x_0)}\times\{0\}\bigr)
\leq
|(z,y)-(x,0)|+\sigma
<\delta.
\end{aligned}
\]
Moreover,
$
|z-x_0|<r+\sigma+4\rho
<r',
\
|y|<\delta/4<h.
$
Consequently,
\[
B_{4\rho}(x,0)\subset D_\delta \subset B_{R/2}(x_0)\times(-\kappa_{0} R,\kappa_{0} R) \ \text{for every} \ x\in B_{r+\sigma}(x_0).
\]
Since $V_y$ satisfies the first equation in
\eqref{2026-xx1}, the interior regularity
estimates for elliptic equations, together
with the Lipschitz continuity of $A(z)$, yield
\begin{equation}\label{eq:W-y-estimate}
\int_{B_{2\rho}(x,0)}
|V_{yy}(z,y)|^2\,dz\,dy
\leq
\frac{C(N,\lambda)}{\rho^2}
\int_{B_{4\rho}(x,0)}
|V_y(z,y)|^2\,dz\,dy.
\end{equation}

On the other hand, applying the
one-dimensional trace inequality ( see \cite[Lemma 9.9]{Brezis}) to
$V_y(z,\cdot)$ on $(-\rho,\rho)$
and integrating over $z\in B_\rho(x)$,
we obtain
\[
\begin{aligned}
\int_{B_\rho(x)}|V_y(z,0)|^2\,dz
\leq{}&
\frac{C(N)}{\rho}
\int_{B_\rho(x)\times(-\rho,\rho)}
|V_y(z,y)|^2\,dz\,dy
+
C(N)\rho
\int_{B_\rho(x)\times(-\rho,\rho)}
|V_{yy}(z,y)|^2\,dz\,dy.
\end{aligned}
\]
Since
$
B_\rho(x)\times(-\rho,\rho)
\subset B_{2\rho}(x,0),
$
combining the above estimate with
\eqref{eq:W-y-estimate} and recalling that
$V_y(\cdot,0)=u(\cdot,\tau)$, we deduce
\begin{equation}\label{eq:local-u}
\int_{B_\rho(x)}|u(z,\tau)|^2\,dz
\leq
\frac{C(N,\lambda)}{\rho}
\int_{B_{4\rho}(x,0)}
|\nabla_{z,y}V(z,y)|^2\,dz\,dy.
\end{equation}

To obtain the estimate on
$B_{r+\sigma}(x_0)$, choose a maximal
$\rho/2$-separated family
$\{x_j\}_{j=1}^{m}$ in
$\overline{B_{r+\sigma}(x_0)}$.
By maximality,
\begin{equation}\label{e2026-11}
B_{r+\sigma}(x_0)
\subset
\bigcup_{j=1}^{m}B_\rho(x_j).
\end{equation}
Moreover, the balls
$\{B_{\rho/4}(x_j)\}_{j=1}^{m}$
are pairwise disjoint.

We claim that the enlarged balls
$B_{4\rho}(x_j,0)$ have uniformly
bounded overlap, namely,
\begin{equation}\label{eq:bounded-overlap}
\sum_{j=1}^m
\mathbf{1}_{B_{4\rho}(x_j,0)}(X)
\leq 17^N
\qquad\text{for all }X\in\mathbb{R}^{N+1}.
\end{equation}
To verify this, fix
$X=(z,y)\in\mathbb{R}^{N+1}$.
If $X\in B_{4\rho}(x_j,0)$,
then $|z-x_j|<4\rho$, and hence
$
B_{\rho/4}(x_j)
\subset B_{17\rho/4}(z).
$
Since these smaller balls are pairwise
disjoint, a comparison of their volumes yields
\[
\mbox{Card}\bigl\{j\in\mathbb{N}^+: X\in B_{4\rho}(x_j,0)
\bigr\}
\leq
\frac{|B_{17\rho/4}(z)|}
{|B_{\rho/4}(x_j)|}
=17^N.
\]
This establishes \eqref{eq:bounded-overlap}.

Combining \eqref{eq:local-u}, \eqref{e2026-11}, and \eqref{eq:bounded-overlap}, we obtain
\[
\begin{aligned}
\|u(\tau)\|_{L^2(B_{r+\sigma}(x_0))}^2
&\leq
\sum_{j=1}^{m}
\int_{B_\rho(x_j)}
|u(z,\tau)|^2\,dz\\
&\leq
\frac{C(N,\lambda)}{\rho}
\sum_{j=1}^{m}
\int_{B_{4\rho}(x_j,0)}
|\nabla_{z,y}V(z,y)|^2\,dz\,dy\\
&=\frac{C(N,\lambda)}{\rho}
\int_{D_\delta}
\left(
\sum_{j=1}^m
\mathbf{1}_{B_{4\rho}(x_j,0)}(z,y)
\right)
|\nabla_{z,y}V(z,y)|^2\,dzdy\\
&\leq
\frac{C(N,\lambda)}{\rho}
\int_{D_\delta}
|\nabla_{z,y}V(z,y)|^2\,dz\,dy.
\end{aligned}
\]
Recalling that $\rho=\delta/16$,
we conclude that
\begin{equation*}
\|u(\tau)\|_{L^2(B_{r+\sigma}(x_0))}
\leq
C(N,\lambda)\delta^{-1/2}
\|V\|_{H^1(D_\delta)}.
\end{equation*}
This completes the proof of \eqref{yu-7-4-1777}.

It follows from (\ref{yu-7-4-1777}) and \eqref{yu-7-4-1} that
\begin{equation}\label{yu-7-5-1}
\begin{aligned}
\|u(\tau)\|_{L^2(B_{r+\sigma}(x_0))}
\leq
C
\|u(\tau)\|_{L^2(B_r(x_0))}^{\alpha}
\|V\|_{L^2(
B_{2r'}(x_0)\times(-2h,2h)
)}^{1-\alpha},
\end{aligned}
\end{equation}
where $C:=C(r,N,\lambda,\alpha)>0$.

To finish the proof, it suffices to bound the term $\|V\|_{L^2(B_{2r'}(x_0)\times(-2h,2h))}$. We will treat
$V_1$ and $V_2$ separately.
On the one hand, by \eqref{yu-6-23-6jia},
Lemma \ref{yu-lemma-6-18-1}, and
Minkowski's inequality, we obtain,
for every $|y|\leq 2h$,
\begin{align}
\|V_1(\cdot,y)\|_{L^2(B_{2r'}(x_0))}
&\leq
C_{3}h(T+R^2)
e^{\frac{C_{1}T}{T-\tau}}
F_T(R)
\int_{\mathbb R}
e^{-q\sqrt{|\mu|}}\,d\mu
\nonumber\\
&=
4C_{3}hq^{-2}
(T+R^2)
e^{\frac{C_{1}T}{T-\tau}}
F_T(R),
\label{yu-7-5-30}
\end{align}
where $ q:=R/C_4-\sqrt2h>0.$
Since $R=C(N,\lambda)r$ is fixed, by \eqref{yu-7-5-30} and integrating over $|y|<2h$, we obtain
\begin{equation}
\|V_1\|_{L^2(
B_{2r'}(x_0)\times(-2h,2h)
)}^2 \leq
C(r,N,\lambda)
(1+T^2)
e^{\frac{2C_1T}{T-\tau}}
F^{2}_T(R).\label{yu-7-5-3}
\end{equation}

On the other hand, by (\ref{yu-7-5-7}) and Parseval's identity we obtain
\begin{equation}\label{yu-7-5-10}
\begin{array}{lll}
\displaystyle{}\|V_2\|^{2}_{L^2(
B_{2r'}(x_0)\times(-2h,2h)
)}
&\leq&\displaystyle{}
\sum_{k=1}^\infty\alpha_k^2e^{-2\mu_k\tau}
\int_{-2h}^{2h}\frac{\sinh^2(\sqrt{\mu_k}y)}
{\mu_k}dy\\
&\leq&\displaystyle{}
C(r,N,\lambda)R^2e^{\frac{2h^2}{\tau}}
\sum_{k=1}^\infty\alpha_k^2\\
&\leq&\displaystyle{}
C(r,N,\lambda)R^2e^{\frac{2h^2}{\tau}}
F^{2}_T(R).
\end{array}
\end{equation}
Therefore, by (\ref{yu-7-5-3}) and (\ref{yu-7-5-10}), we conclude that
\begin{equation*}
\|V\|_{L^2(
B_{2r'}(x_0)\times(-2h,2h)
)}\leq C(r,N,\lambda)
\left[
(1+T^2)
e^{\frac{2C_1T}{T-\tau}}
+
e^{\frac{2h^2}{\tau}}
\right]^{\frac{1}{2}}F_{T}(R).
\end{equation*}
This, together with (\ref{yu-7-5-1}) and $h=r$, means that
\[
\begin{aligned}
\|u(\cdot,\tau)\|_{L^2(B_{r+\sigma}(x_0))}
\leq C
\left[
(1+T^2)
e^{\frac{2C_1T}{T-\tau}}
+
e^{\frac{2r^2}{\tau}}
\right]^{\frac{1-\alpha}{2}}
\|u(\cdot,\tau)\|_{L^2(B_r(x_0))}^{\alpha}
F_T(R)^{1-\alpha},
\end{aligned}
\]
with $C=C(N,r,\alpha,\lambda)>0$. This proves \eqref{eq-local-ball} and completes the proof.
\end{proof}

\begin{Lemma}\label{yu-LEMMA-7-10-6}
For each $\theta\in(0,1)$ and $r>0$, there exist
positive constants
$C_7=C_7(\theta,r,N,\lambda)$
and
$\sigma=\sigma(\theta,r,N,\lambda)
\in(0,\min\{1,r/2\})$
such that, for every $0<T<1$ and every
$u_0\in H^1(\mathbb R^N)$,
the corresponding solution $u$ of
\eqref{yu-6-24-1} satisfies
\begin{equation}\label{eq-global-ball}
\begin{aligned}
\|u(\cdot,T)\|_{L^2(\mathbb R^N)}
\leq{}&
C_7e^{r^{2}\theta/T}
\|u(\cdot,T)\|_{L^2(\omega^\sigma)}^{1-\theta}
\left(
\sup_{s\in[0,2T]}
\|u(\cdot,s)\|_{H^1(\mathbb R^N)}
\right)^\theta,
\end{aligned}
\end{equation}
where $\omega^\sigma$ is the observation set
defined in \eqref{2026-3a}.
\end{Lemma}

\begin{proof}
Set $ a_\sigma=2(r+\sigma/2)/\sqrt N$ and $x_k=a_\sigma k,\ k\in\mathbb Z^N.$
We first establish three geometric properties
of the lattice balls. First, we show that
the enlarged balls $B_{r+\sigma}(x_k)$,
$k\in\mathbb Z^N$, cover the whole space
$\mathbb R^N$. 

For any $x=(x_1,\ldots,x_N)\in\mathbb R^N$,
there exists
$k=(k_1,\ldots,k_N)\in\mathbb Z^N$
such that
$
|x_j-a_\sigma k_j|
\leq a_\sigma/2, \ j=1,\ldots,N.
$
Consequently,
\[
|x-x_k|
\leq\frac{\sqrt N a_\sigma}{2}
<r+\sigma.
\]
It follows that
\begin{equation}\label{eqq}
\mathbb R^N
=\bigcup_{k\in\mathbb Z^N}
B_{r+\sigma}(x_k).
\end{equation}

Second, we prove that the
families $\{B_S(x_k)\}_{k\in\mathbb Z^N}$
have uniformly bounded overlap for every
fixed $S>0$. In particular, this applies
to the families $\{B_r(x_k)\}_k$ and
$\{B_{2R}(x_k)\}_k$, which will be used
to derive the global interpolation inequality
from the local estimate. Let $S>0$ and fix
$x\in\mathbb R^N$.
If $x\in B_S(x_k)$, then
$
|x_j-a_\sigma k_j|<S,
\ j=1,\ldots,N,
$
which implies
\[
\frac{x_j-S}{a_\sigma}
<k_j<
\frac{x_j+S}{a_\sigma}.
\]
Since an interval of length $2S/a_\sigma$
contains at most $1+2S/a_\sigma$
integers and $
a_\sigma\geq2r/\sqrt N
$, we obtain
\begin{equation}\label{eq-cover-overlap}
\sup_{x\in\mathbb R^N}
\sum_{k\in\mathbb Z^N}
\mathbf1_{B_S(x_k)}(x)
\leq
\left(1+\frac{2S}{a_\sigma}\right)^N
\leq
\left(1+\frac{S\sqrt N}{r}\right)^N.
\end{equation}
Thus, for fixed $r,S,N$, the overlap
bound is independent of $\sigma$.

Finally, we show
that the complement of the observation set
$\omega^\sigma$ contains the open balls
$B_{\sigma/2}(p_k)$, whose centers tend
to infinity as $|k|\to\infty$.
For each $k\in\mathbb Z^N$, define
\[
p_k=a_\sigma
\left(k+\left(\frac12,\ldots,\frac12\right)\right).
\]
For every $j\in\mathbb Z^N$,
\[
\begin{aligned}
|p_k-x_j|^2
&=a_\sigma^2
\sum_{i=1}^N
\left(k_i+\frac12-j_i\right)^2
\geq\frac{Na_\sigma^2}{4}.
\end{aligned}
\]
Equality is attained when $j=k$.
Hence
$
\operatorname{dist}
(p_k,a_\sigma\mathbb Z^N)
=r+\sigma/2.
$
For any $z\in B_{\sigma/2}(p_k)$
and $j\in\mathbb Z^N$,
\[
\begin{aligned}
|z-x_j|
&\geq|p_k-x_j|-|z-p_k|>r.
\end{aligned}
\]
Therefore,
$
B_{\sigma/2}(p_k)
\subset
\mathbb R^N\setminus\omega^\sigma.
$

We now prove the interpolation
inequality \eqref{eq-global-ball}.

By Lemma \ref{lemma-2A2}, with
$\alpha=1-\theta$, for every
$k\in\mathbb Z^N$ we have
\begin{equation}\label{eqqas}
\|u(\cdot,\tau)\|_{L^2(B_{r+\sigma}(x_k))}\leq
C_6
\left[
(1+T^2)e^{\frac{2C_1T}{T-\tau}}
+e^{\frac{2r^2}{\tau}}
\right]^{\theta/2}
\|u(\cdot,\tau)\|_{L^2(B_r(x_k))}^{1-\theta}
\left(
\sup_{0\leq s\leq T}
\|u(\cdot,s)\|_{H^1(B_R(x_k))}
\right)^\theta,
\end{equation}
where $R=C(N,\lambda)r$ and
$C_6=C_6(\theta,r,N,\lambda)$.
For $0<T\leq 2$ and $0<\tau\leq T/2$,
we have
\[
(1+T^2)e^{\frac{2C_1T}{T-\tau}}
+e^{\frac{2r^2}{\tau}}
\leq
(5e^{4C_1}+1)e^{\frac{2r^2}{\tau}}.
\]
This, together with \eqref{eqqas}, \eqref{eqq} and
H\"older's inequality for sequences, implies that
\begin{equation}\label{3.44444}
\begin{aligned}
\|u(\cdot,\tau)\|_{L^2(\mathbb R^N)}^2
\leq{}&
C(r,N,\lambda,\theta)e^{2r^2\theta/\tau}
\left(
\sum_{k\in\mathbb Z^N}
\|u(\cdot,\tau)\|_{L^2(B_r(x_k))}^2
\right)^{1-\theta}\\
&\times
\left(
\sum_{k\in\mathbb Z^N}
\sup_{s\in[0,T]}
\|u(\cdot,s)\|_{H^1(B_R(x_k))}^2
\right)^\theta.
\end{aligned}
\end{equation}
By \eqref{eq-cover-overlap} with
$S=r$ and Tonelli's theorem,
\begin{equation}\label{3.4X12}
\begin{aligned}
&\sum_{k\in\mathbb Z^N}
\|u(\cdot,\tau)\|_{L^2(B_r(x_k))}^2=
\int_{\omega^\sigma}
|u(x,\tau)|^2
\sum_{k\in\mathbb Z^N}
\mathbf1_{B_r(x_k)}(x)\,dx\leq
(1+\sqrt N)^N\|u(\cdot,\tau)\|_{L^2(\omega^\sigma)}^2.
\end{aligned}
\end{equation}

Next, we show that
\begin{equation}\label{yu-4-9-4}
\sum_{k\in\mathbb Z^N}\sup_{s\in[0,T]}\|u(\cdot,s)\|^{2}_{H^1(B_{R}(x_k))}\leq C(1+T)\sup_{s\in[0,T]}\|u(\cdot,s)\|^2_{H^1(\mathbb{R}^{N})},
\end{equation}
where $C:=C(r,N,\lambda)>0$.

Indeed, choose a cutoff function
$\psi\in C_0^\infty(\mathbb R^N;[0,1])$
such that $\operatorname{supp}\psi\Subset B_{2R}(0)$
and $\psi\equiv 1$ in $B_R(0)$.
For each $k\in\mathbb Z^N$, set
$\psi_k(x):=\psi(x-x_k)$ and
$w_k(x,t):=\psi_k(x)u(x,t)$.
Then
$\operatorname{supp}\psi_k\Subset B_{2R}(x_k)$
and $\psi_k\equiv 1$ in $B_R(x_k)$.
A direct computation shows that $w_k$
satisfies
\begin{equation}\label{yu-4-9-2}
\begin{cases}
(w_k)_t+\mathcal{A}w_k
=-2\nabla\psi_k\cdot(A(x)\nabla u)
+u\mathcal{A}\psi_k
&\text{in } B_{2R}(x_k)\times(0,T),\\
w_k=0
&\text{on } \partial B_{2R}(x_k)\times(0,T),\\
w_k(\cdot,0)=\psi_k u_0
&\text{in } B_{2R}(x_k).
\end{cases}
\end{equation}
By the standard energy estimate for (\ref{yu-4-9-2}), one can easily check that
\begin{equation}\label{eq-local-energy}
\begin{aligned}
&\sup_{s\in[0,T]}
\|u(\cdot,s)\|_{H^1(B_R(x_k))}^2\leq \sup_{s\in[0,T]}
\|w_k(\cdot,s)\|_{H^1(B_{2R}(x_k))}^2\\
\leq&
C(r,N,\lambda)
\left[
\|u_0\|_{H^1(B_{2R}(x_k))}^2
+\int_0^T
\|u(\cdot,s)\|_{H^1(B_{2R}(x_k))}^2\,ds
\right].
\end{aligned}
\end{equation}

We now sum \eqref{eq-local-energy}
over $k\in\mathbb Z^N$.
By \eqref{eq-cover-overlap}
with $S=2R$,
\[
\begin{aligned}
\sup_{x\in\mathbb R^N}
\sum_{k\in\mathbb Z^N}
\mathbf1_{B_{2R}(x_k)}(x)
\leq
\left(1+\frac{2R\sqrt N}{r}\right)^N.
\end{aligned}
\]
Consequently, for every $s\in[0,T]$, we have
\[
\begin{aligned}
\sum_{k\in\mathbb Z^N}
\|u(s)\|_{H^1(B_{2R}(x_k))}^2\leq
C(r,N,\lambda)
\|u(s)\|_{H^1(\mathbb R^N)}^2.
\end{aligned}
\]
Summing \eqref{eq-local-energy} over \(k\in\mathbb Z^N\) and applying Tonelli's theorem, we obtain \eqref{yu-4-9-4}.

Finally, by
\eqref{yu-4-9-4}, \eqref{3.4X12} and 
\eqref{3.44444}, we obtain
\[
\begin{aligned}
\|u(\cdot,\tau)\|_{L^2(\mathbb R^N)}^2
\leq{}
C(\theta,r,N,\lambda)
e^{2r^2\theta/\tau}
\|u(\cdot,\tau)\|_{L^2(\omega^\sigma)}^{2(1-\theta)}
\left(
\sup_{s\in[0,T]}
\|u(\cdot,s)\|_{H^1(\mathbb R^N)}^2
\right)^\theta.
\end{aligned}
\]
Replacing $\tau$ and $T$ in the above
inequality by $T$ and $2T$, respectively,
we obtain, for every $T\in(0,1)$,
\[
\begin{aligned}
\|u(\cdot,T)\|_{L^2(\mathbb R^N)}^2
\leq{}&
C(\theta,r,N,\lambda)e^{2r^2\theta/T}
\|u(\cdot,T)\|_{L^2(\omega^\sigma)}^{2(1-\theta)}
\left(
\sup_{s\in[0,2T]}
\|u(\cdot,s)\|_{H^1(\mathbb R^N)}^2
\right)^\theta.
\end{aligned}
\]
Taking square roots yields
\eqref{eq-global-ball}, which completes
the proof.
\end{proof}

\textbf{Proof of Theorem \ref{yu-theorem-7-10-6}.}
By standard energy estimates of solutions to \eqref{yu-6-24-1}, we have
\begin{equation}\label{yu-7-12-2}
	\|u(\cdot,t)\|_{H^1(\mathbb R^N)}\leq \frac{Ce^{Ct}}{\sqrt{t}}\|u_0\|_{L^2(\mathbb R^N)}
	\end{equation}
with $C=C(N,\lambda)>0$, for each $t\in(0,6T]$. Moreover, if $u_0\in H^1(\mathbb R^N)$, then
 \begin{equation}\label{yu-7-12-3}
 	\|u(\cdot,t)\|_{H^1(\mathbb R^N)}\leq Ce^{Ct}\|u_0\|_{H^1(\mathbb R^N)}\;\;\mbox{for each}\;\;
	t\in[0,6T],
 \end{equation}
where $C=C(N,\lambda)>0$.

We consider the following equation
\begin{equation*}
\begin{cases}
	v_t-\mbox{div}(A(x)\nabla v)=0&\mbox{in}\;\mathbb R^N\times(0,4T),\\
	v(\cdot,0)=u(\cdot,\frac{T}{2})&\mbox{in}\;\mathbb R^N.
\end{cases}
\end{equation*}
	It is clear that $v(\cdot,t)=u(\cdot,t+T/2)$ when $t\in[0,4T]$. Moreover, by (\ref{yu-7-12-2}) we have
$u(\cdot,T/2)\in H^1(\mathbb R^N)$,
which means that $v\in C([0,4T];H^1(\mathbb R^N))$.
For any \(\varepsilon\in(0,1)\), applying Lemma \(\ref{yu-LEMMA-7-10-6}\) with \(u_0\) and \(T\) replaced by \(u(\cdot,T/2)\) and \(T/2\), respectively, yields
\begin{equation*}
	\left\|v\left(\cdot,\frac{T}{2}\right)\right\|_{L^2(\mathbb R^N)}\leq C_7e^{2r^{2}\varepsilon/T} \left\|v\left(\cdot,\frac{T}{2}\right)\right\|^{1-\varepsilon}_{L^2(\omega^{\sigma})}
\left(\sup_{s\in[0,T]}\|v(\cdot,s)\|_{H^1(\mathbb R^N)}\right)^{\varepsilon}.
\end{equation*}
This, along with (\ref{yu-7-12-3}), gives that
\begin{equation*}
\left\|v\left(\cdot,\frac{T}{2}\right)\right\|_{L^2(\mathbb R^N)}\leq C_7e^{2r^{2}\varepsilon/T}e^{C\varepsilon T} \left\|v\left(\cdot,\frac{T}{2}\right)\right\|^{1-\varepsilon}_{L^2(\omega^{\sigma})}
\|v(\cdot,0)\|^{\varepsilon}_{H^1(\mathbb R^N)}.
\end{equation*}
Thus,
\begin{equation*}
\left\|u\left(\cdot,T\right)\right\|_{L^2(\mathbb R^N)}\leq  C_{8}e^{C_{9}\varepsilon(T+\frac{1}{T})} \left\|u\left(\cdot,T\right)\right\|^{1-\varepsilon}_{L^2(\omega^{\sigma})}
\left\|u\left(\cdot,\frac{T}{2}\right)\right\|^{\varepsilon}_{H^1(\mathbb R^N)}
\end{equation*}
where $C_{8}=C_{8}(r,N,\lambda,\varepsilon)>0$ and $C_{9}=C_{9}(r,N,\lambda)>0$.
Combining this with (\ref{yu-7-12-2}) and using \(T\in(0,1)\), we obtain
\begin{equation}\label{yu-8-20}
\left\|u\left(\cdot,T\right)\right\|_{L^2(\mathbb R^N)}\leq  C_{8}e^{\frac{C_{10}\varepsilon}{T}} \left\|u\left(\cdot,T\right)\right\|^{1-\varepsilon}_{L^2(\omega^{\sigma})}
\left\|u_{0}\right\|^{\varepsilon}_{L^2(\mathbb R^N)}
\end{equation}
where $C_{10}=C_{10}(r,N,\lambda)>0$.
This completes the proof of
Theorem \ref{yu-theorem-7-10-6}.
\qed

\section{Proof of Theorem \ref{jiudu4}}\label{finalproof}

\textbf{Proof of Theorem \ref{jiudu4}.}
Fix $\varepsilon\in(0,1)$ and $T\in(0,1)$. We apply
Theorem~\ref{yu-theorem-7-10-6} with an auxiliary parameter
$\theta\in(0,1/2)$, to be specified at the end of the proof, and set
$\alpha=\theta/(1-\theta)\in(0,1)$.
By \eqref{yu-8-20} and Young's inequality, 
for any $0\le t_1<t_2\le T$ and every $\gamma>0$,
\begin{equation}\label{2019-7-9}
\begin{aligned}
\|u(\cdot,t_2)\|_{L^2(\mathbb R^N)}^2
\leq {}&\gamma\|u(\cdot,t_1)\|_{L^2(\mathbb R^N)}^2+C_\alpha\gamma^{-\alpha}
\exp\!\left(\frac{2C_{10}\alpha}{t_2-t_1}\right)
\|u(\cdot,t_2)\|_{L^2(\omega^{\sigma_\theta})}^2.
\end{aligned}
\end{equation}
Here $C_{10}=C_{10}(r,N,\lambda)>0$ is the constant from
\eqref{yu-8-20}, while $C_\alpha$ may
depend on $\alpha,r,N,\lambda$ and may change from line to line.
According to Proposition 2.1 in \cite{Phung-Wang-2013}, for each $\kappa > 1$ the sequence $\{l_m\}_{m\in \mathbb{N}^{+}}$, give by
$$
 l_m:=\frac{T}{\kappa^{m-1}},\ \text{for each}\ m\in \mathbb{N}^{+}.
$$
Then $l_m\downarrow0$ and
\[
 l_m-l_{m+1}=\frac{T(\kappa-1)}{\kappa^m}.
\]
Fix $m\ge1$ and $t\in[l_{m+1},l_m]$. Applying
\eqref{2019-7-9} with $t_1=l_{m+2}$ and $t_2=t$, we obtain
\begin{equation*}
\begin{aligned}
\|u(\cdot,t)\|_{L^2(\mathbb R^N)}^2
\leq {}&\gamma\|u(\cdot,l_{m+2})\|_{L^2(\mathbb R^N)}^2+C_\alpha\gamma^{-\alpha}
\exp\!\left(\frac{2C_{10}\alpha}{t-l_{m+2}}\right)
\|u(\cdot,t)\|_{L^2(\omega^{\sigma_\theta})}^2.
\end{aligned}
\end{equation*}
By a standard energy estimate, we have that
$$\|u(\cdot,l_m)\|_{L^2(\mathbb R^N)}^2\le
\|u(\cdot,t)\|_{L^2(\mathbb R^N)}^2.$$
Since $t-l_{m+2}\ge l_{m+1}-l_{m+2}$, integration over
$t\in(l_{m+1},l_m)$ gives
\begin{equation}\label{3.2525253}
\begin{aligned}
(l_m-l_{m+1})\|u(\cdot,l_m)\|_{L^2(\mathbb R^N)}^2
\leq& \gamma(l_m-l_{m+1})
\|u(\cdot,l_{m+2})\|_{L^2(\mathbb R^N)}^2\\
&\quad+C_\alpha\gamma^{-\alpha}
\exp\!\left(\frac{2C_{10}\alpha}{l_{m+1}-l_{m+2}}\right)
\int_{l_{m+1}}^{l_m}
\|u(\cdot,t)\|_{L^2(\omega^{\sigma_\theta})}^2\,dt.
\end{aligned}
\end{equation}
Therefore, dividing \eqref{3.2525253} by $l_m-l_{m+1}$ yields
\begin{equation}\label{3.2525254}
\begin{aligned}
\|u(\cdot,l_m)\|_{L^2(\mathbb R^N)}^2
\leq {}&\gamma\|u(\cdot,l_{m+2})\|_{L^2(\mathbb R^N)}^2
+\frac{C_\alpha}{e\kappa C_{10}\alpha}\gamma^{-\alpha}
\exp\!\left(\frac{3C_{10}\alpha\kappa^{m+1}}{T(\kappa-1)}\right)
\int_{l_{m+1}}^{l_m}
\|u(\cdot,t)\|_{L^2(\omega^{\sigma_\theta})}^2\,dt.
\end{aligned}
\end{equation}
Set
$
 d:=3C_{10}\alpha/[T\kappa(\kappa-1)].
$
Multiplying \eqref{3.2525254} by
$\gamma^\alpha e^{-d\kappa^{m+2}}$ and choosing
$\gamma=e^{-d\kappa^{m+2}}$, we obtain
\begin{equation}\label{3.25252555}
\begin{aligned}
e^{-(1+\alpha)d\kappa^{m+2}}
\|u(\cdot,l_m)\|_{L^2(\mathbb R^N)}^2
-e^{-(2+\alpha)d\kappa^{m+2}}
\|u(\cdot,l_{m+2})\|_{L^2(\mathbb R^N)}^2
\leq\frac{C_\alpha}{e\kappa C_{10}\alpha}
\int_{l_{m+1}}^{l_m}
\|u(\cdot,t)\|_{L^2(\omega^{\sigma_\theta})}^2\,dt.
\end{aligned}
\end{equation}
Now choose
$
 \kappa:=\sqrt{(\alpha+2)/(\alpha+1)}>1.
$
Then $(1+\alpha)\kappa^2=2+\alpha$, and \eqref{3.25252555}
becomes
\[
\begin{aligned}
e^{-(2+\alpha)d\kappa^m}
\|u(\cdot,l_m)\|_{L^2(\mathbb R^N)}^2
-e^{-(2+\alpha)d\kappa^{m+2}}
\|u(\cdot,l_{m+2})\|_{L^2(\mathbb R^N)}^2
\leq\frac{C_\alpha}{e\kappa C_{10}\alpha}
\int_{l_{m+1}}^{l_m}
\|u(\cdot,t)\|_{L^2(\omega^{\sigma_\theta})}^2\,dt.
\end{aligned}
\]
Summing the preceding inequality over \(m\ge1\) and using the telescoping property, we obtain
\begin{eqnarray*}
&&e^{-(2+\alpha)d\kappa^{2}}\|u(\cdot,T)\|^{2}_{L^{2}(\mathbb{R}^{N})}
\leq e^{-(2+\alpha)d\kappa^{2}}\|u(\cdot,l_{2})\|^{2}_{L^{2}(\mathbb{R}^{N})}\\
&\leq&\sum_{m=1}^{+\infty}\left(e^{-(2+\alpha)d\kappa^{m}}\|u(\cdot,l_{m})\|^{2}_{L^{2}(\mathbb{R}^{N})}
-e^{-(2+\alpha)d\kappa^{m+2}}\|u(\cdot,l_{m+2})\|^{2}_{L^{2}(\mathbb{R}^{N})}\right)\\
&\leq&\frac{C_\alpha}{e\kappa C_{10}\alpha}
\int_{0}^{T}\|u(\cdot,t)\|^{2}_{L^{2}(\omega^{\sigma_\theta})}\mathrm{d}t.
 \end{eqnarray*}
Since $0<\alpha<1$, we have
$\kappa\ge\sqrt{3/2}$ and $\alpha+2\le3$. Hence
\[
\frac{3C_{10}\alpha(\alpha+2)}{T(\kappa-1)}
\leq\frac{\alpha K}{T},\qquad
K:=\frac{9C_{10}}{\sqrt{3/2}-1},
\]
where $K$ is independent of $\alpha$ and $T$.
Consequently,
\[
\|u(\cdot,T)\|_{L^2(\mathbb R^N)}^2
\leq C_\alpha e^{\alpha K/T}
\int_0^T\|u(\cdot,t)\|_{L^2(\omega^{\sigma_\theta})}^2\,dt.
\]
Finally, for the prescribed $\varepsilon\in(0,1)$, choose
\[
\alpha:=\min\left\{\frac12,\frac{\varepsilon}{2K}\right\},
\qquad \theta:=\frac{\alpha}{1+\alpha}\in(0,1/2),
\]
and take $\sigma=\sigma_\theta$ from
Theorem~\ref{yu-theorem-7-10-6}.
Since $\alpha K\le\varepsilon$, it follows that
\[
\int_{\mathbb R^N}|u(x,T)|^2\,dx
\leq C(\varepsilon,r,N,\lambda)e^{\varepsilon/T}
\int_0^T\int_{\omega^\sigma}|u(x,t)|^2\,dx\,dt,
\]
where $C(\varepsilon,r,N,\lambda)>0$ is independent of $T$.
This completes the proof of Theorem~\ref{jiudu4}.\qed


\end{document}